\documentclass{cimart}

\usepackage[markup=underlined]{changes}

\usepackage[shortlabels]{enumitem}
\usepackage[all]{xy}

\newcommand{\be }{\begin{equation}}
\newcommand{\ee }{\end{equation}}

\newcommand{\mf}{\mathfrak m}
\newcommand{\Nat}{\mathbb N}
\newcommand{\huaR}{\mathcal{R}}
\newcommand{\huaO}{\mathcal{O}}
\newcommand{\huaN}{\mathcal{N}}

\newcommand{\Id}{\mathrm{Id}}
\newcommand{\Hom}{\mathrm{Hom}}
\newcommand{\gl}{\mathfrak {gl}}
\newcommand{\ad}{\mathrm{ad}}

\VOLUME{32}
\ISSUE{2}
\YEAR{2024}
\firstpage{127}
\DOI{https://doi.org/10.46298/cm.12656}

\title{
    Cohomology and deformations of left-symmetric Rinehart algebras
    }

\author{
    Abdelkader  Ben Hassine, Taoufik Chtioui, Mohamed Elhamdadi and Sami Mabrouk
    }

\authorinfo[
    A. Ben Hassine]{
    Department of Mathematics, University of Bisha, Saudi Arabia}{%
    Benhassine@ub.edu.sa
    }

\authorinfo[
    T. Chtioui]{
    Faculty of Sciences, Gabes University, Tunisia}{%
    chtioui.taoufik@yahoo.fr
    }

\authorinfo[
    M. Elhamdadi]{
    Department of Mathematics, University of South Florida, U.S.A.}{%
    emohamed@math.usf.edu
    }

    \authorinfo[
    S. Mabrouk]{
    Faculty of Sciences, University of Gafsa, Tunisia}{%
    sami.mabrouk@fsgf.u-gafsa.tn, mabrouksami00@yahoo.fr
    }

\abstract{
We introduce a notion of left-symmetric Rinehart algebras, which is a generalization of the notion of left-symmetric algebras. The left multiplication gives rise to a representation of the corresponding sub-adjacent Lie-Rinehart algebra. We construct left-symmetric Rinehart algebras from $\mathcal O$-operators on Lie-Rinehart algebras. 
We extensively investigate representations of left-symmetric Rinehart algebras.  Moreover, we construct a graded Lie algebra on the space of multi-derivations whose Maurer–Cartan elements characterize left-symmetric Rinehart algebras and study deformations of left-symmetric Rinehart algebras, which are controlled by the second cohomology class in the deformation cohomology. We also give  the relationships between $\mathcal O$-operators and Nijenhuis operators on left-symmetric Rinehart algebras.
}

\keywords{
    left-symmetric Rinehart algebra,  representation, graded Lie algebra, Maurer-Cartan element,  cohomology, deformation, Nijenhuis operator.
    }

\msc{
    17E05 (primary); 53D17, 17B70, 14B12, 06B15 (secondary)
    }

\begin{document}

\section{Introduction}

Left-symmetric algebras are algebras for which the associator 
\[(x,y,z):=(x \cdot y) \cdot z-x \cdot (y \cdot z)\]
satisfies the identity $(x,y,z)=(y,x,z)$.  These algebras appeared as early as 1896 in the work of Cayley \cite{Cayley} as  
rooted tree algebras.
In the 1960s, they also arose from the study of several topics in geometry and algebra, such as convex
homogenous cones \cite{Vinberg}, affine manifolds and affine structures on Lie groups \cite{Koszul,Matsushima} and deformations of
associative algebras \cite{Gerstenhaber}. 
In 2006, Burde \cite{Burde}
wrote an interesting survey 
showing the importance of left-symmetric algebras 
in many areas, such
as vector fields, rooted tree algebras, 
vertex algebras, operad theory, deformation
complexes of algebras, convex homogeneous cones, affine manifolds and left-invariant affine structures on
Lie groups \cite{Burde}.

 Left symmetric algebras are the underlying algebraic structures of 
non-abelian phase spaces of Lie algebras
\cite{Kuper,non-abelian}, 
leading to a
bialgebra theory of left-symmetric algebras \cite{Left-symmetric-bialgebras}. They can also be seen 
as the algebraic structures
 behind the classical Yang-Baxter equations.  Precisely, they provide a
construction of solutions of the classical Yang-Baxter equations in
certain semidirect product Lie algebra structures (that is, over the
double spaces) induced by left-symmetric algebras
\cite{Kuper2,Bai:CYBE}.

The notion of Lie-Rinehart algebras was introduced by J. Herz in \cite{Herz1953} and further
developed in \cite{Palais, Rinehart2}. A a notion of (Poincar\'e) duality for this class of algebras was introduced in \cite{Huebschmann,Huebschmann1}.
 Lie-Rinehart structures have been the subject of extensive studies, 
 in relations to symplectic geometry, Poisson structures, Lie groupoids and algebroids   and other kinds of quantizations (see \cite{H2004,Krahmer&Rovi,LuiShengBai1,LuiShengBai2,LuiShengBai3,Mackenzie}).  For further details and a history of the notion of
Lie-Rinehart algebra, we refer the reader to \cite{H2004}. Lie-Rinehart algebras have been investigated
furthermore in \cite{Bkouche,Casas,Chen&Liu&Zhong,Chen&Liu&Liu}.

A left-symmetric algebroid is a geometric generalization of a left-symmetric algebra. See
\cite{LuiShengBai1,LuiShengBai2,LuiShengBai3} for more details and applications. The notion of a Nijenhuis operator on a left-symmetric
algebroid was introduced in \cite{LuiShengBai2}, which could generate a trivial deformation. More details on deformations of left-symmetric algebras can be found in  \cite{WangBaiLiuSheng}.

In this paper, we introduce a notion of  left-symmetric Rinehart algebras, which is a generalization of a left-symmetric algebra and an algebraic version of left symmetric algebroids.   The following diagram shows  how left-symmetric Rinehart algebras fit in relation to Lie algebras, left-symmetric algebras and Lie-Rinehart algebras.

\[
    \xymatrix{
        \mbox{Lie algebra} \ar[d] \ar[rrr]^{\mbox{\small generalization}}
        &&& \mbox{Lie-Rinehart algebra} \ar[d] \\
        \mbox{Left-symmetric algebra} \ar@<-1ex>[u] \ar[rrr]^{\mbox{\small generalization}}
        &&& \mbox{Left-symmetric Rinehart} \ar@<-1ex>[u]
    }
\]

The paper is organized as follows. In Section 2, we recall some definitions concerning  left-symmetric algebras
 and Lie-Rinehart algebra. In Section 3, we introduce the  notion of left-symmetric Rinehart algebra and give some of its properties.  As in the case of a left-symmetric algebras, one can obtain the sub-adjacent Lie-Rinehart algebra from a left-symmetric Rinehart algebra by using the commutator.
 The left multiplication gives rise to a representation of the sub-adjacent Lie-Rinehart algebra. We construct left-symmetric Rinehart algebras using $\mathcal O$-operators.  Section 4 is devoted to the study of representations and cohomology  of left-symmetric Rinehart algebra. In Section 5, we construct a graded Lie algebra whose Maurer-Cartan elements are left-symmetric
Rinehart algebras which give rise to a coboundary operator. Section 6 is devoted to
   introduce the deformation cohomology associated to a left-symmetric Rinehart algebra, which 
  controls the deformations. In
  Section 7, we  introduce the notion of a Nijenhuis operator,
 which could generate a trivial deformation. In addition, we investigate some connection 
 between $\mathcal O$-operators and Nijenhuis operators.

 Throughout this paper all vector spaces are over a field $\mathbb{K}$ of characteristic  zero. 

\section{Preliminaries}

In this section, we briefly recall some   basics  of left-symmetric algebras and  Lie-Rinehart algebras \cite{Burde}.

\begin{definition}  A  left-symmetric algebra   is a vector space $L$ endowed with a linear map $\cdot:L\otimes L\longrightarrow L$
such that for any $x,y,z\in L$,
$$(x,y,z)=(y,x,z),\;\;{\rm or}\;\;{\rm
equivalently,}\;\;(x\cdot y)\cdot z-x\cdot(y\cdot z)=(y\cdot x)\cdot
z-y\cdot(x\cdot z),$$
where the associator
$(x,y,z):=(x\cdot y)\cdot z-x\cdot(y\cdot z)$.
\end{definition}
 Let $\ad^L$ (resp. $\ad^R$)
be the left multiplication operator (resp. right multiplication operator) on $L$ that is, i.e. $\ad^L(x)y=x\cdot
y$ (resp. $\ad^R(x)y=y\cdot
x$), for any $x,y\in L$. The following lemma is given in \cite{Burde}.
\begin{lemma} \label{lem:sub-ad} Let $(L,\cdot)$ be a left-symmetric algebra. The commutator
$ [x,y]=x\cdot y-y\cdot x$  defines a Lie algebra $L$,
which is called the  sub-adjacent Lie algebra  of $L$.  The algebra $L$ is also
called a  compatible left-symmetric algebra  on the Lie
algebra $L$. Furthermore, the map $ad^L:L\rightarrow
\gl(L)$ with $x \mapsto  L_x$ gives a representation of the Lie
algebra 
$(L,[\cdot,\cdot])$.
\end{lemma}

\begin{definition}Let $(L,\cdot)$ be a left-symmetric algebra and $M$ a vector space. A  representation of $L$ on $M$ consists of a pair $(\rho,\mu)$,
 where $\rho:L\longrightarrow \gl(M)$ is a representation of the  sub-adjacent Lie algebra $L$ on $M$ and $\mu:L\longrightarrow \gl(M)$ is a linear map satisfying:
\begin{eqnarray}\label{rho-mu}
 \rho(x)\circ\mu(y)-\mu(y)\circ\rho(x)=\mu(x\cdot y)-\mu(y)\circ\mu(x),\quad \forall~x,y\in L.
\end{eqnarray}
  \end{definition}
The map $\rho$ is called a left representation and 
$\mu$ is a right representation.
Usually, we denote a representation by $(M;\rho,\mu)$.  Then $(L;\ad^L,\ad^R)$ is a representation of $(L,\cdot)$ which is called adjoint representation.

The cohomology complex for a left-symmetric algebra $(L,\cdot)$ with a representation $(M;\rho,\mu)$ is given as follows. The set of
$(n+1)$-cochains is given by
\begin{equation}
 C^{n+1}(L,M)=\Hom(\wedge^n L\otimes L,M),\quad
 \forall n\geq 0.
\end{equation}
For all $\omega\in C^n(L,M)$, the coboundary operator $\delta:C^n(L,M)\longrightarrow C^{n+1}(L,M)$ is given by
\begin{eqnarray*}
\delta \omega(x_1,x_2,\dots,x_{n+1})
&=&\sum_{i=1}^n(-1)^{i+1}\rho(x_i)\omega(x_1,\dots,\widehat{x_i},\dots,x_{n+1})\\
&&\quad+\sum_{i=1}^n(-1)^{i+1}\mu(x_{n+1})\omega(x_1,\dots,\widehat{x_i},\dots,x_n,x_i)\\
&&\quad-\sum_{i=1}^n(-1)^{i+1}\omega(x_1,\dots,\widehat{x_i},\dots,x_n,x_i\cdot x_{n+1})\\
&&\quad+\sum_{1\leq i<j\leq n}(-1)^{i+j}\omega([x_i,x_j],x_1,\dots,\widehat{x_i},\dots,\widehat{x_j},\dots,x_{n+1}).
\end{eqnarray*}

We then have the following lemma whose proof comes from a direct computation using identity~(\ref{rho-mu}).
\begin{lemma}[See \cite{Burde2}]
The map $\delta$ satisfies $\delta^2=0$.
\end{lemma}

\begin{definition} \label{Lie-rinehart}
A Lie-Rinehart algebra $L$ over an associative commutative algebra $A$ is a Lie algebra
over $\mathbb K$ with an $A$-module structure and a  linear map $\rho: L \rightarrow Der(A)$, such that the following conditions hold:
 \begin{enumerate}
\item For all $a\in A$ and $ x,y  \in L$
$$\rho([x,y])=\rho(x)\rho(y)-\rho(x)\rho(y)\quad \text{and}\quad \rho(ax) = a\rho(x).$$ 
\item The compatibility condition:
\begin{equation}\label{compatibility-Lie-rinhe}
[x, ay] = \rho(x)ay + a[x, y],~~~\forall a\in A, x,y \in L.
\end{equation}
\end{enumerate}
\end{definition}

Let $(L, A,[\cdot,\cdot]_L, \rho)$ and $(L', A',[\cdot,\cdot]_{L'}, \rho')$ be two Lie-Rinehart algebras,  then a Lie-Rinehart algebra homomorphism is defined as a pair of maps $(g, f)$, where
the maps $f:L\rightarrow L'$ and $g: A \rightarrow A'$ are two  algebra homomorphisms such that:
\begin{itemize}
\item[(1)] $f(ax) = g(a)f(x)$ for all $x \in L$ and $ a \in A,$
\item[(2)] $g(\rho(x)a) = \rho'(f(x))g(a)$ for all $x\in L$ and $ a \in A.$
\end{itemize}

Now, we recall the definition  of  module over a Lie-Rinehart algebra (for more details see \cite{Chemla}).

\begin{definition}\label{rep-Lie rinehart}
Let $M$ be an $A$-module. Then  $M$ is a  module over a
Lie-Rinehart algebra $(L,A,[\cdot,\cdot],\rho)$  if there exists  a map $\theta : L\otimes M \rightarrow M$  such that:
\begin{enumerate}
 \item $\theta$ is a representation of the Lie algebra
$(L, [\cdot,\cdot])$ on $M$.
\item $\theta(ax, m) = a\theta (x, m)$ for all $a \in A, x \in  L, m \in  M.$
\item $\theta(x, am) = a\theta(x, m)+ \rho(x)am$ for all $x\in L, a \in A, m \in M.$
\end{enumerate}
\end{definition}
We have the following lemma giving a characterization of of the $\theta$ which are representations.

\begin{lemma}\label{rep to produit semidirect}
The map $\theta$ is representation if and only if $L\oplus M$ is  Lie-Rinehart algebra over $A$, where 
$[\cdot,\cdot]_{L\oplus M}$ and $\theta_{L\oplus M}$ are given by
\begin{eqnarray*}
[x_1+m_1,x_2+m_2]_{L\oplus M}&=&[x_1, x_2]+\rho(x_1)m_2-\rho(x_2)m_1, \\
\theta_{L\oplus M}(x_1+m_1)&=&\theta(x_1)   
\end{eqnarray*}
for all $x_1,x_2\in L$ and $  m_1,m_2\in M$.
\end{lemma}
%
\section{Some basic properties of a left-symmetric Rinehart algebras}

In this section, we introduce a notion of left-symmetric Rinehart  algebras  illustrated by some examples. As in the case of a left-symmetric
algebra, we obtain the sub-adjacent Lie-Rinehart  algebra from a left-symmetric Rinehart algebra using the
commutator. In addition, we construct left-symmetric Rinehart algebras using $\mathcal O$-operators.
\begin{definition}\label{defi:left-symmetric-Rinehart}

A  left-symmetric Rinehart algebra is a quadruple $(L,A,\cdot,\ell)$ where $(L,\cdot)$ is a left-symmetric algebra, $A$ is an associative commutative algebra  and $\ell: L \rightarrow Der(A)$ a linear map  such  that the following conditions hold:
 \begin{enumerate}
 \item $L$ is an $A$-module.
\item For all $a\in A$ and $ x,y  \in L$
$$\ell(x\cdot y-y\cdot x)=\ell(x)\ell(y)-\ell(x)\ell(y)\quad \text{and}\quad \ell(ax) = a\ell(x).$$ 
\item The compatibility conditions: for all $a\in A$ and $ x,y \in L$
\begin{align}
x\cdot(ay) =& \ell(x)ay + a(x\cdot y),\label{compatibility-left-symmetric-Rinehart1}\\
(ax)\cdot y =& a(x\cdot y).\label{compatibility-left-symmetric-Rinehart2}
\end{align}
\end{enumerate}
\end{definition}

\begin{example}
It is clear that any left-symmetric  algebra is a left-symmetric Rinehart algebra.
\end{example}
\begin{example}
A Novikov Poisson algebra  is a left-symmetric Rinehart algebra 
(see  \cite{Xu2}).
\end{example}

\begin{example}
 Let $(L, A, \cdot,\ell)$ be a left-symmetric Rinehart algebra and let $L\oplus A$
 be the direct sum of $L$ and $A$.
Then $(L\oplus A, A, \cdot_{L\oplus A},\ell_{L\oplus A})$ is a left-symmetric Rinehart algebra, where the $\cdot_{L\oplus A}$ is defined by the following expression, for all $x_1, x_2\in L$,
$a_1, a_2\in A$;
\begin{equation*}
(x_1+ a_1)\cdot_{L\oplus A} (x_2+ a_2) = x_1\cdot x_2 + \ell(x_1)(a_2);
\end{equation*}
and $\ell_{L\oplus A}:L\oplus A \rightarrow Der(A)$ is defined by $\ell_{L\oplus A}(a_1+ x_1) = \ell(x_1)$.  Indeed, it obvious that $(L\oplus A,\cdot_{L\oplus A})$ is a left-symmetric algebra, $\ell_{L\oplus A}$ is a representation of left-symmetric algebra $L\oplus A$
 and $\ell_{L\oplus A}\in Der(A)$.\\
 By direct calculation, we have $\ell_{L\oplus A}(b(x_1+ a_1))=b\ell_{L\oplus A}(x_1+ a_1)$ for all $b,a_1\in A $ and $x_1\in L$.
On the other hand, letting $x_1,x_2\in L$ and $b,a_1,a_2\in A$, we have
\begin{align*}
(x_1+ a_1)\cdot_{L\oplus A} b(x_2+ a_2) =&(a_1+ x_1)\cdot_{L\oplus A} (bx_2+ ba_2)\\
=&x_1\cdot (bx_2) + \ell(x_1)(ba_2)\\
=&b(x_1\cdot x_2) + \ell(x_1)b(x_2)+\ell(x_1)(b)a_2+b\ell(x_1)(a_2)\\
=&b\big(x_1\cdot x_2+\ell(x_1)(a_2)\big)+\ell(x_1)b(x_2)+\ell(x_1)(b)a_2\\
=&b\big((x_1+ a_1)\cdot_{L\oplus A} (x_2+ a_2)\big)+\ell_{L\oplus A}(x_1+ a_1)b(x_2+ a_2).
\end{align*}
Moreover,
\begin{align*}
b(x_1+ a_1)\cdot_{L\oplus A} (x_2+ a_2) =&(bx_1+ ba_1)\cdot_{L\oplus A} (x_2+ a_2)\\
=&(bx_1)\cdot x_2 + \ell(bx_1)(a_2)\\
=&b(x_1\cdot x_2) + b\ell(x_1)(a_2)\\
=&b\big(x_1\cdot x_2 + \ell(x_1)a_2\big)\\
=&b\big((x_1+ a_1)\cdot_{L\oplus A} (x_2+ a_2)\big).
\end{align*}
\end{example}


Now we have the following theorem.

\begin{theorem}\label{theorem:sub-adjacent}
  Let $(L,A,\cdot,\ell)$ be a left-symmetric Rinehart algebra.
Then, $(L, A,[\cdot,\cdot],\ell)$   
 is a Lie-Rinehart algebra,  denoted by
$L^C$, called the  sub-adjacent Lie-Rinehart algebra of
 $(L,A,\cdot,\ell)$.
\end{theorem}

\begin{proof}
Since $(L,\cdot)$ is a left-symmetric  algebra,  we
have that $(L,[\cdot,\cdot])$ is a Lie algebra. For any
$a\in A$, by direct computations, we have
\begin{eqnarray*}
  [x,ay]&=&x\cdot(ay)-(ay)\cdot x=a(x\cdot y)+\ell(x)ay-a(y\cdot x)\\
  &=&a[x,y]+\ell(x)(a)y,
\end{eqnarray*}
which implies that $(L,A,[\cdot,\cdot],\ell)$ is a Lie-Rinehart algebra.

To see that the linear map $\ell:L\longrightarrow Der(A)$ is a
representation, we only need to show that
$\ell_{[x,y]}=[\ell_x,\ell_y]_{Der(A)}$, which follows directly from the fact
that $(L,\cdot)$ is a left-symmetric  algebra. This ends the proof.\end{proof}

\begin{definition}
Let $(L_1, A_1,\cdot_1,\ell_1)$ and $(L_2,A_2,\cdot_2,\ell_2)$ be  two left-symmetric
Rinehart algebras. A  homomorphism  of left-symmetric Rinehart algebras is a pair of two algebra homomorphisms  $(f,g)$  where $f:L_1\longrightarrow L_2$ and $g:A_1\longrightarrow A_2$  such that:
$$f(ax)=g(a)f(x), \:~g(\ell_1(x)a)=\ell_2(f(x))g(a), \quad\forall x,y\in L_1, a\in A_1. $$
\end{definition}
The following proposition is immediate.
\begin{proposition}
Let $(f,g)$ be a homomorphism of
left-symmetric Rinehart algebras from  $(L_1, A_1,\cdot_1,\ell_1)$ to $(L_2,A_2,\cdot_2,\ell_2)$. Then $(f,g)$ is also a Lie-Rinehart algebra
homomorphism of the corresponding sub-adjacent Lie-Rinehart algebras.
\end{proposition}

Now we give the definition of an $\mathcal{O}$-operator.

\begin{definition}
  Let $(L,A,[\cdot,\cdot],\rho)$ be a Lie-Rinehart algebra and $\theta:L\longrightarrow End(M)$ be a representation over $M$. A linear map
  $T:M\longrightarrow L$ is called an  $\mathcal{O}$-operator if for all $u,v\in M$ and $a\in A$ we have
  \begin{equation}\label{O-operator1}
  T(au)=aT(u),
  \end{equation}
 \begin{equation}\label{O-operator2}
  [T(u),T(v)]=T\big(\theta(T(u))(v)-\theta(T(v))(u)\big).
 \end{equation}
\end{definition}

\begin{remark}
  Consider the semidirect product Lie-Rinehart algebra 
  \[(L\ltimes_{\theta} M,A,[\cdot,\cdot]_{L\ltimes_{\theta} M},\rho_{L\ltimes_{\theta} M}),\]
   where $\rho_{L\ltimes_{\theta} M}(x+u):=\rho(x)(u)$ and the bracket $[\cdot,\cdot]_{L\ltimes_{\theta} M}$ is given by
  $$
  [x+u,y+v]_{L\ltimes_{\theta} M}=[x,y]+\theta(x)(v)-\theta(y)(u).
  $$
  Any $\huaO$-operator $T:M\longrightarrow L$ gives a Nijenhuis operator  $\tilde{T}=\left(\begin{array}{cc}
    0&T\\0&0
  \end{array}\right)$  on the Lie-Rinehart algebra $L\ltimes_\theta
 M$. More precisely, we have
  {\small$$
  [\tilde{T}(x+u),\tilde{T}(y+v)]_{L\ltimes_{\theta} M}=\tilde{T}\Big([\tilde{T}(x+u),y+v]_{L\ltimes_{\theta} M}+[x+u,\tilde{T}(y+v)]_{L\ltimes_{\theta} M}-\tilde{T}[x+u,y+v]_{L\ltimes_{\theta} M}\Big).
  $$}
 Fore more details on Nijenhuis operators and their applications the reader should consult~\cite{Dorfman}.
\end{remark}

Let $T:M\longrightarrow L$ be an $\huaO$-operator. Define the multiplication $\cdot_T$ on $M$ by
$$
u\cdot_T  v=\theta(T(u))(v), \: \forall u,v\in M.
$$
We then have the following proposition.

\begin{proposition}
  With the above notations, $(M,A, \cdot_T,\ell_T=\ell\circ T)$ is a left-symmetric Rinehart algebra, and the map $T$ is Lie-Rinehart algebra homomorphism from 
  $(M,[\cdot,\cdot])$ to $(L,[\cdot,\cdot])$.
\end{proposition}
\begin{proof} It is easy to see that  $(M,\cdot_T)$ is a left-symmetric
algebra. For any $a\in A$, using Definition \ref{defi:left-symmetric-Rinehart} and equation \eqref{O-operator1} we have
\begin{equation*}
\ell_M(au)=\ell(T(au))=a\ell(T(u))=a\ell_M(u),
\end{equation*}
Similarly, using Definition \ref{rep-Lie rinehart} we obtain
\begin{eqnarray*}
  (au)\cdot_T v&=&\theta(T(au))(v)=\theta aT(u))(v)=a\theta(T(u))(v),\\
  u\cdot_T(av)&=&\theta(T(u))(av)=a\theta(T(u))(v)+\ell\circ T(u)(a)v.
\end{eqnarray*}
Thus, $(M,A, \cdot_T,\ell_M)$ is a left-symmetric Rinehart algebra.  Let 
 $[\cdot,\cdot]$  be the sub-adjacent Lie bracket on $M$. Then
we have
$$
T[u,v]=T(u\cdot_T v-v\cdot_T
u)=T(\theta(T(u))(v)-\theta(T(v))(u))=[T(u),T(v)].
$$
So $T$ is a homomorphism of  Lie algebras.\end{proof}

\section{Representations  of left-symmetric Rinehart algebras}
In this section, we develop the notion of representations of a left-symmetric Rinehart algebra and give a cohomology theory with coefficients in a representation.

\begin{definition}
Let $(L,A,\cdot,\ell)$ be a left-symmetric Rinehart algebra and $M$ be an $A$-module. A  representation of $A$ on $M$ consists of a pair
$(\rho,\mu)$, where $\rho$ is a representation
of the sub-adjacent Lie-Rinehart algebra $(L,A,[\cdot,\cdot]^C,\ell)$ and $\mu:L  \rightarrow End(M)$ is a linear
map, such that for all $x,y\in L$ and $\ m\in M$, we have
\begin{eqnarray}\label{representation condition 2}
&&\mu(ax)m=a\mu(x)m=\mu(x)(am)\nonumber\\
 &&\rho(x)\mu(y)-\mu(y)\rho(x)=\mu(x\cdot y)-\mu(y)\mu(x).
\end{eqnarray}
We will denote this representation by $(M;\rho,\mu)$.
\end{definition}

For a left-symmetric Rinehart algebra $(L,A,\cdot,\ell)$ and a representation $(M;\rho,\mu)$, the following proposition gives a construction of a left-symmetric Rinehart algebra called semidirect product and denoted by $L\ltimes_{\rho,\mu} M$.
\begin{proposition}\label{produit-semidirect}
Let $(L,A,\cdot,\ell)$ be a left-symmetric Rinehart algebra and $(M;\rho,\mu)$ a representation. Then, $(L\oplus M,A,\cdot_{L\oplus M},\ell_{L\oplus M})$ is a left-symmetric Rinehart algebra, where
$\cdot_{L\oplus M}$ and $\ell_{L\oplus M}$ are given by
\begin{eqnarray}
(x_1+m_1)\cdot_{L\oplus M}(x_2+m_2)&=&x_1\cdot x_2+\rho(x_1)m_2+\mu(x_2)m_1,\label{eq:mulsemi}\\
\ell_{L\oplus M}(x_1+m_1)&=&\ell(x_1),\label{eq:mulsemi1}
\end{eqnarray}
for all $x_1,x_2\in L$ and $\ m_1,m_2\in M$.
\end{proposition}

\begin{proof} Let $(M;\rho,\mu)$ be a representation. It is straightforward to  see that
$(L\oplus M,A,\cdot_{L\oplus M})$ is a left-symmetric algebra. For any $x_1,x_2\in L$ and $m_1,m_2\in M$, we have
\begin{eqnarray*}
(x_1+m_1)\cdot_{L\oplus M} (a(x_2+m_2))&=&x_1\cdot (ax_2)+\rho(x_1)am_2+\mu(ax_2)m_1\\
&=& a(x_1\cdot x_2)+\ell(x_1)(ax_2)+a\rho(x_1) m_2 \\
&&\qquad\qquad\:+\ell(x_1)(a m_2)+a\mu(x_2)m_1\\
&=& a((x_1+m_1)\cdot_{L\oplus M}  (x_2+m_2))+\ell_{L\oplus M}(x_1)(a)(x_2+m_2).
\end{eqnarray*}On the other hand, we have
 \begin{eqnarray*}
(a(x_1+m_1))\cdot_{L\oplus M} (x_2+m_2)&=&(ax_1)\cdot x_2+\rho(ax_1)m_2+\mu(x_2)(am_1)\\
&=& a((x_1+m_1)\cdot_{L\oplus M}(x_2+m_2)).
\end{eqnarray*}
Therefore, $(L\oplus M,A,\cdot_{L\oplus M},\ell_{L\oplus M})$ is a left-symmetric Rinehart algebra.
\end{proof}

Let $(L,A,\cdot,\ell)$ be a left-symmetric Rinehart algebra  and
$(M;\rho,\mu)$ be a representation.\\
Let $\rho^*:L\otimes M^* \rightarrow M^*$  and  $\mu^*:M^*\otimes L \longrightarrow M^*$ be defined by
\begin{align*}
\langle \rho^*(x)\xi,m\rangle=\ell(x)\langle \xi,m\rangle -\langle \rho(x)m,\xi \rangle\ \ \textrm{and}\ \
\langle \mu^*(x)\xi,m\rangle=-\langle\xi, \mu(x)m\rangle,
\end{align*}
where $M^*=Hom_A(M,A)$.
Then, we  have the following proposition.

\begin{proposition}\label{pro:representation}
With the above notations, we obtain that
\begin{itemize}
\item[$\rm(i)$] $(M,\rho-\mu)$ is a representation of the sub-adjacent Lie-Rinehart algebra $(L,A,[\cdot,\cdot],\ell)$.
\item[$\rm(ii)$] $(M^*,\rho^*-\mu^*,-\mu^*)$ is a representation be a left-symmetric Rinehart algebra $(L,A,\cdot,\ell)$.
\end{itemize}
\end{proposition}
\begin{proof} Since $(M;\rho,\mu)$ is a representation of the left-symmetric
algebra $(L,A,\cdot,\ell)$, using Proposition \ref{produit-semidirect} we have that $(L\oplus M,A,\cdot_{L\oplus M},\ell_{L\oplus M})$ is a left-symmetric Rinehart algebra. Consider its sub-adjacent Lie-Rinehart algebra $(L\oplus M,A,\cdot_{L\oplus M},[\cdot,\cdot]_{L\oplus M},\ell_{L\oplus M})$. We have
\begin{eqnarray*}
[(x_1+m_1),(x_2+m_2)]_{L\oplus M}&=&(x_1+m_1)\cdot_{L\oplus M}(x_2+m_2)-(x_2+m_2)\cdot_{L\oplus M}(x_1+m_1)\\
&=&x_1\cdot x_2+\rho(x_1)m_2+\mu(x_2)m_1 \\
&&-x_2\cdot x_1-\rho(x_2)m_1-\mu(x_1)m_2\\
&=&[x_1,x_2]^C+(\rho-\mu)(x_1)(m_2)+(\rho-\mu)(x_2)(m_1).
\end{eqnarray*}
From Lemma \ref{rep to produit semidirect} we deduce that $(M,\rho-\mu)$ is a representation of Lie-Rinehart algebra $L$ on $M$. This finishes the proof of (i).

For item (ii), it is clear that $\rho^*-\mu^*$ is just the dual representation of the representation $\theta=\rho-\mu$ of the sub-adjacent Lie-Rinehart algebra of $L$.
We can directly check that $-\mu^*(ax)\xi=-a\mu^*(x)\xi=-\mu^*(x)(a\xi)$. For any $x,y\in L$, $\xi\in M^*$ and $m\in M$ we have
\small{
\begin{eqnarray*}
&&\hspace{-1cm}-\langle (\rho^*-\mu^*)(x)\mu^*(y)\xi,m\rangle+\langle\mu^*(y)(\rho^*-\mu^*)(x)\xi,m\rangle \\
&&=
-\langle \rho^*(x)\mu^*(y)\xi,m\rangle+\langle \mu^*(x)\mu^*(y)\xi,m\rangle
+\langle \mu^*(y)\rho^*(x)\xi,m\rangle-\langle \mu^*(y)\mu^*(x)\xi,m\rangle\\
&&=\ell(x)\langle \xi,\mu(y)m\rangle-\langle \xi,\mu(y)\rho(x)m\rangle
+\langle \xi,\mu(y)\mu(x)m\rangle\\
&&\qquad-\ell(x)\langle \xi, \mu(y)m\rangle+\langle \xi, \rho(x)\mu(y)m\rangle-\langle \xi,\mu(x)\mu(y)m\rangle\\
&&=\langle \xi,\mu(x\cdot y)m\rangle-\langle \xi, \mu(y)\mu(x)m\rangle
+\langle \xi,\mu(y)\mu(x)m\rangle-\langle \xi,\mu(x)\mu(y)m\rangle\\
&&=\langle(-\mu^*(x\cdot y)-\mu^*(y)\mu^*(x))\xi,m\rangle.
\end{eqnarray*}
}
Therefore $(M^*,\rho^*-\mu^*,-\mu^*)$ is a representation of $L$.
\end{proof}

\begin{corollary} With the above notations, we have
  \begin{itemize}
\item[$\rm(i)$] The left-symmetric Rinehart algebras $L\ltimes_{\rho,\mu} M$ and
$L\ltimes_{\rho-\mu,0}M$  have the same sub-adjacent
Lie-Rinehart algebra $L\ltimes_{\rho-\mu}M.$
\item[$\rm(ii)$] The left-symmetric Rinehart algebras
$L\ltimes_{\rho^*,0} M^*$ and $L\ltimes_{\rho^*-\mu^*,-\mu^*} M^*$
have the same sub-adjacent Lie-Rinehart algebra $L\ltimes_{\rho^*}M^*$.
\end{itemize}
\end{corollary}

Let $(M;\rho,\mu)$ be a   representation of a left-symmetric
Rinehart algebra $(L,A,\cdot,\ell)$. In general, $(M^*,\rho^*,\mu^*) $ is not
a representation. But we have the following proposition.
\begin{proposition}\label{dual representation condition}
Let  $(L,A,\cdot,\ell)$  be a left-symmetric Rinehart algebra and
$(M;\rho,\mu)$ be a representation. Then the following conditions
are equivalent:
\begin{itemize}
\item[$\rm(1)$]$(M;\rho-\mu,-\mu)$ is a representation of $(L,A,\cdot,\ell)$.
\item[$\rm(2)$]$(M^*,\rho^*,\mu^*)$ is a representation of $(L,A,\cdot,\ell)$.
\item[$\rm(3)$]$\mu(x)\mu(y)=\mu(y)\mu(x)$ for all $x,y\in L.$
\end{itemize}
\end{proposition}



\section{The Matsushima-Nijenhuis bracket for left-symmetric Rinehart algebras}
 
In this section we construct a graded Lie algebra whose Maurer-Cartan elements are left-symmetric Rinehart algebras which give rise to a coboundary operator.

Let    $(L,A,\cdot,\ell)$ be a left-symmetric Rinehart algebras.
A multiderivation of degree $n$ is a multilinear map $P \in \Hom(\Lambda^{n}L\otimes L,L)$ such that, for every $a \in A$, $x_i\in L$ and $i \in \{1,2, \dotsc, n+1\}$, we have  
\begin{align}
&\label{condition-cochain}
P(x_1,\dots,ax_i,\dots, x_n,x_{n+1}) =aP(x_1, \dots, x_i,\dots, x_n,x_{n+1}), \\
&P(x_1,\cdots,x_n,ax_{n+1})=aP(x_1,\cdots,x_n,x_{n+1})+\Xi_{P} (x_1,\cdots, x_{n})(a)x_{n+1},  
 \end{align}
 where $\Xi_P: L^{\otimes n} \to Der(A)$ is called the symbol map. 
 The space of all multiderivations of degree $n$ will be denoted by $\mathfrak D^n(L)$. Set 
 $\mathfrak{D}^*(L)=\displaystyle\oplus_{n\geq -1} \mathfrak D^{n}(L)$ with   $\mathfrak D^{-1}(L)=L$, the space of multiderivations on $L$. 
 
 A permutation $\sigma\in\mathbb S_n$ is called an $(i,n-i)$-unshuffle if $\sigma(1)<\cdots<\sigma(i)$ and $\sigma(i+~1)<~\cdots<\sigma(n)$. If $i=0$ and $i=n$, we assume $\sigma=Id$. The set of all $(i,n-i)$-unshuffles will be denoted by $\mathbb S_{(i,n-i)}$. The notion of an $(i_1,\cdots,i_k)$-unshuffle and the set $\mathbb S_{(i_1,\cdots,i_k)}$ are defined similarly.

Let $P \in  \mathfrak D^{m}(L)$ and $Q\in \mathfrak D^{n}(L)$. We define the Matsushima--Nijenhuis bracket 
$[\cdot, \cdot]_{MN} :~\mathfrak D^{m}(L)\times \mathfrak D^{n}(L)\to \mathfrak D^{m+n}(L) $ by  \[  [P, Q]_{MN} = P\diamond  Q- (-1)^{mn} Q\diamond  P,\]
where 

\begin{eqnarray*}
&&P\diamond Q(x_1,x_2,\cdots,x_{m+n+1})\\
&&=\sum_{\sigma\in \mathbb{S}_{(m,1,n-1)}}(-1)^{\sigma}P(Q(x_{\sigma(1)},\cdots,x_{\sigma(m+1)}),x_{\sigma(m+2)},\cdots,x_{\sigma(m+n)},x_{m+n+1})\\
&&\quad+(-1)^{mn}\sum_{\sigma\in \mathbb{S}_{(n,m)}}(-1)^{\sigma}P(x_{\sigma(1)},\cdots,x_{\sigma(n)},Q(x_{\sigma(n+1)},x_{\sigma(n+2)},\cdots,x_{\sigma(m+n)},x_{m+n+1})).
\end{eqnarray*}

\begin{theorem}\label{MN bracket}
With the above notations, we have 
\begin{itemize}
    \item [(i)] The pair $(\mathfrak D^*(L), [\cdot,\cdot]_{MN})$ is a graded Lie algebra. 
    \item [(ii)]  There is a one-to-one correspondence between the set of   Maurer-Cartan elements of the graded Lie algebra $(\mathfrak D^*(L),[\cdot,\cdot]_{MN})$ and 
    left-symmetric Rinehart algebra structures on $L$.  
\end{itemize}
\end{theorem}

\begin{proof}
(i)\ We begin by check that the Matsushima-Nijenhuis bracket is well defined. For $P\in \mathfrak D^m(L)$ and $Q\in  \mathfrak D^n(L)$, by a direct calculation, we have
{\small \begin{eqnarray*}
&&\hspace{-0.8cm}[P,Q]_{MN}(ax_1,x_2,\cdots,x_{m+n+1})\\
&&\hspace{-0.3cm}=aP\diamond Q(x_1,x_2,\cdots,x_{m+n+1})-(-1)^{mn}aQ\diamond P(x_1,x_2,\cdots,x_{m+n+1})\\
&&+\sum_{\sigma\in \mathbb{S}_{(m-1,1,n-1)}}(-1)^{\sigma}\Xi_{Q}(x_{\sigma(2)},\cdots,x_{\sigma(m+1)})(a)P(x_1,x_{\sigma(m+2)},\cdots,x_{\sigma(m+n)},x_{m+n+1})\\
&&+(-1)^{mn}\sum_{\sigma\in \mathbb{S}_{(n-1,1,m-1)}}(-1)^{\sigma}\Xi_{P}(x_{\sigma(2)},\cdots,x_{\sigma(n+1)})(a)Q(x_1,x_{\sigma(n+2)},\cdots,x_{\sigma(m+n)},x_{m+n+1})\\
&&-(-1)^{mn}\sum_{\sigma\in \mathbb{S}_{(n-1,1,m-1)}}(-1)^{\sigma}\Xi_{P}(x_{\sigma(2)},\cdots,x_{\sigma(n+1)})(a)Q(x_1,x_{\sigma(n+2)},\cdots,x_{\sigma(m+n)},x_{m+n+1})\\
&&-\sum_{\sigma\in \mathbb{S}_{(m-1,1,n-1)}}(-1)^{\sigma}\Xi_{Q}(x_{\sigma(2)},\cdots,x_{\sigma(m+1)})(a)P(x_1,x_{\sigma(m+2)},\cdots,x_{\sigma(m+n)},x_{m+n+1})\\
&&\hspace{-0.3cm}=a[P,Q]_{MN}(x_1,x_2,\cdots,x_{m+n+1}),
\end{eqnarray*}}
which implies that
\begin{eqnarray*}
[P,Q]_{MN}(ax_1,x_2,\cdots,x_{m+n+1})=a[P,Q]_{MN}(x_1,x_2,\cdots,x_{m+n+1}).
\end{eqnarray*}
It is straightforward to check that $[P,Q]_{MN}$ is skew-symmetric with respect to its first $m+n$ arguments. Thus $[P,Q]_{MN}$ is $A$-linear with respect to its first $m+n$ arguments.

On the other hand, following a straightforward  calculation, we have
\begin{eqnarray*}
[P,Q]_{MN}(x_1,x_2,\cdots,ax_{m+n+1})
&=&a[P,Q]_{MN}(x_1,x_2,\cdots,x_{m+n+1})\\
&&+\Xi_{[P,Q]_{MN}}(x_1,x_2,\cdots,x_{m+n})(a)x_{m+n+1},
\end{eqnarray*}
where the symbol map $\Xi_{[P,Q]_{MN}}$ is given by
\begin{eqnarray*}
&&\hspace{-2cm}\Xi_{[P,Q]_{MN}}(x_1,x_2,\cdots,x_{m+n})(a) \\
=&&\sum_{\sigma\in \mathbb{S}_{(m,1,n-1)}}(-1)^{\sigma}\Xi _P(Q(x_{\sigma(1)},\cdots,x_{\sigma(m+1)}),x_{\sigma(m+2)},\cdots,x_{\sigma(m+n)}))(a) \\
&+&\sum_{\sigma\in \mathbb{S}_{(n,1,m-1)}}(-1)^{\sigma}\Xi_Q(P(x_{\sigma(1)},\cdots,x_{\sigma(n+1)}),x_{\sigma(n+2)},\cdots,x_{\sigma(m+n)})(a) \\
&+&(-1)^{mn}\sum_{\sigma\in \mathbb{S}_{(m,n)}}(-1)^{\sigma}\Xi _P(x_{\sigma(1)},\cdots,x_{\sigma(n)})(\Xi_Q (x_{\sigma(n+1)},\cdots,x_{\sigma(n+m)}))(a)\\
&+&\sum_{\sigma\in \mathbb{S}_{(m,n)}}(-1)^{\sigma}\Xi _Q(x_{\sigma(1)},\cdots,x_{\sigma(m)})(\Xi_P (x_{\sigma(m+1)},\cdots,x_{\sigma(m+n)}))(a).
\end{eqnarray*}
Thus $[P,Q]_{MN}\in\mathfrak D^{m+n}(L)$.

It was shown in \cite{ChaLiv} that the Matsushima-Nijenhuis bracket provides a graded Lie algebra structure on the graded vector space  $\oplus_{n\geq 1} Hom(\Lambda^{n-1}L\otimes L,L)$.  Therefore,   $(\mathfrak D^*(L),[\cdot,\cdot]_{MN})$ is a graded Lie algebra.

(ii)\ Let  $\pi\in \mathfrak D^1(L)$, we have
\begin{eqnarray*}
  \pi(ax_1,x_2)=a\pi(x_1,x_2),\quad  \pi(x_1,ax_2)=a\pi(x_1,x_2)+\Xi_\pi(x_1)(a)x_2,\quad \forall~x_1,x_2\in L.
\end{eqnarray*}
In addition, we can easily check that   
\begin{eqnarray*}
[\pi,\pi]_{MN}(x_1,x_2,x_3)&=&2(\pi(\pi(x_1,x_2),x_3)-\pi(\pi(x_2,x_1),x_3)\\
&&-\pi(x_1,\pi(x_2,x_3))+\pi(x_2,(x_1,x_3))).
\end{eqnarray*}
Thus $(L,A, \pi,\Xi_\pi)$ is a left-symmetric Rinehart algebra if and only if $[\pi,\pi]_{MN}=0$.
\end{proof}

\begin{remark}
 The cohomology of left-symmetric algebras first appeared in the unpublished
paper of Y. Matsushima. Then A. Nijenhuis constructed a graded Lie bracket, which
produces the cohomology theory for left-symmetric algebras. Thus the aforementioned graded
Lie bracket is usually called the Matsushima–Nijenhuis bracket.
\end{remark}

Let $(L,A,\pi,\ell)$ be a left-symmetric Rinehart algebra. According to Theorem \ref{MN bracket}, we have $[\pi,\pi]_{MN}=0$. Using  the graded
Jacobi identity, we get a coboundary operator $\delta: \mathfrak D^{n-1}(L) \to \mathfrak D^n(L)$, by putting 
\begin{align}
    & \delta(P)=(-1)^{n-1}[\pi,P]_{MN},\quad \forall P \in \mathfrak D^{n-1}(L). 
\end{align}
By straightforward computation, we obtain

\begin{proposition}
For any $ P \in \mathfrak D^{n-1}(L)$, we have 
{\small 
\begin{eqnarray}
\delta P(x_1,x_2,\cdots,x_{n+1})
 &=&\sum_{i=1}^{n}(-1)^{i+1}\pi(x_i,P(x_1,x_2,\cdots,\hat{x_i},\cdots,x_{n+1})) \nonumber \\
\label{eq:deformation complex} &&+\sum_{i=1}^{n}(-1)^{i+1}\pi(P(x_1,x_2,\cdots,\hat{x_i},\cdots,x_n,x_i), x_{n+1})\\
\nonumber &&-\sum_{i=1}^{n}(-1)^{i+1}P(x_1,x_2,\cdots,\hat{x_i},\cdots,x_n,\pi(x_i, x_{n+1}))\\
\nonumber &&+\sum_{1\leq i<j\leq {n}}(-1)^{i+j}P(\pi(x_i,x_j)-\pi(x_j,x_i),x_1,\cdots,\hat{x_i},\cdots,\hat{x_j},\cdots,x_{n+1})
\end{eqnarray}}
for all $\ x_i\in L,i=1,2\cdots,n+1$ and $\Xi_{\delta P}$
is given by
{\small
\begin{eqnarray}
\nonumber\Xi_P(x_1,x_2,\cdots,x_n)
&=& \sum_{i=1}^{n}(-1)^{i+1}[\Xi_\pi(x_i),\Xi_{P}(x_1,x_2,\cdots,\hat{x_i},\cdots,x_{n})]\\
\nonumber&&+\sum_{1\leq i<j\leq n}(-1)^{i+j}\Xi_{P}(\pi(x_i,x_j)-\pi(x_j,x_i),x_1,\cdots,\hat{x_i},\cdots,\hat{x_j},\cdots,x_n)\\
&&+\sum_{i=1}^{n}(-1)^{i+1}\Xi_\pi(P(x_1,x_2,\cdots,\hat{x_i},\cdots,x_n,x_i)).\label{eq:simbol}
\end{eqnarray}
}
\end{proposition}

\begin{definition}
The cochain complex $(\mathfrak D^*(L)=\oplus _{n\geq 0}\mathfrak D^n(L),\delta)$ is called
the  deformation complex of the left-symmetric Rinehart algebra $L$. The corresponding $k$-th cohomology group, which we denote by $H^k(L)$, is called the $k$-th deformation cohomology group.
\end{definition}

\section{Deformation of left-symmetric Rinehart algebra}

We investigate in this section a deformation theory of  left-symmetric Rinehart algebras.  But first let us introduce some notation. For a left-symmetric Rinehart algebra  $(L,A,\mathfrak m,\ell)$ we will denote the left-symmetric multiplication $``\cdot "$ by $\mathfrak m$ in the sequel of the paper.
 Let $\mathbb{K}[[t]]$ be the formal power series ring in one variable $t$ and coefficients in  $\mathbb{K}$.  Let $L[[t]] $ be the set of formal power series whose coefficients are elements of $L$ (note that  $L[[t]] $ is
obtained by extending the coefficients domain of $\mathbb{K}[[t]]$ from $\mathbb{K}$ to $L$).  Thus,  $L[[t]] $ is a $\mathbb{K}[[t]]$-module.

\subsection{Formal deformations}

\begin{definition}
  A deformation  of a left-symmetric Rinehart algebra  $(L,A,\mathfrak m,\ell)$ is a $\mathbb{K}[[t]]$-bilinear map
$$\mathfrak{m}_t:L[[t]] \otimes L[[t]]\rightarrow L[[t]]$$
which is given by $\mathfrak{m}_t(x,y)=\displaystyle\sum_{i\geq0} t^i \mathfrak m_i(x,y), $ where $\mathfrak m_0=\mathfrak{m}$ and  the $\mathfrak m_i\in \mathfrak D^1(L)$ satisfy the condition $[\mathfrak m_t,\mathfrak m_t]_{MN}=0$.
\end{definition}
Note that $\mathfrak{m}_t$ is a $1$-degree multiderivation of $A$ with symbol $\Xi_{\mathfrak m_t}: L \to Der(A)$ given  by
$$
    \Xi_{\mathfrak m_t}=\sum_{i\geq0} t^i\Xi_{\mathfrak m_i}.
$$
Moreover, since $[\mathfrak m_t,\mathfrak m_t]_{MN}=0$, it corresponds to a left symmetric Lie Rinehart algebra  structure. In particular, it yields a $t$-parameterized family of products
$
\mathfrak m_t:L\otimes L\to L 
$
and a family of  maps
$
\ell_t:L \rightarrow Der(A),
$
which satisfy the following identities for all $x,y\in L$:
\begin{align*}
\mf_t(x,y)=&x\cdot y+\displaystyle\sum_{i \geq 1}t^i\mf_i(x,y),\\
  \ell_t(x)=&\ell(x)+ \displaystyle\sum_{i \geq 1} t^i \Xi_{\mf_i} (x).
\end{align*}
 The $t$-parametrized family $(L,A,\mf_t, \ell_t)$ is called a $1$-parameter formal deformation of $(L,A,\mf,\ell)$ generated by $\mf_1,\cdots,\mf_m \in \mathfrak D^1(L)$. 

Let $(L,A,\mf_t, \ell_t)$ be a deformation of $\mathfrak{m}$. Then, for all $a\in A, x,y,z \in L$ 
\begin{align}
   &\mathfrak m_t(\mathfrak m_t(x,y),z)-\mathfrak m_t(x,\mathfrak m_t(y,z))=\mathfrak m_t(\mathfrak m_t(y, x),z)-\mathfrak m_t(y,\mathfrak m_t(x, z)).\label{condition} \\
   & \mf_t(ax,y)=a\mf_t(x,y),\label{com2}\\
&\mf_t(x,ay)=a\mf_t(x,y)+\ell_t(x)ay\label{com3}
\end{align}
The identities  \eqref{com2}--\eqref{com3}, mean that $\mf_i\in \mathfrak D^1(L)$. 
Comparing the coefficients of $t^n$ for $n\geq 0$ in equation \eqref{condition}, we get the following:
\begin{equation}\label{DefEqu}
\sum_{i+j=n}\mathfrak m_i(\mathfrak m_j(x,y),z)-\mathfrak m_i(x,\mathfrak m_j(y,z))-\mathfrak m_i(\mathfrak m_j(y, x),z)+ \mathfrak m_i(y,\mathfrak m_j(x, z))=0.
\end{equation}
For $n=1$, equation \eqref{DefEqu} implies
\begin{align*}
&\mathfrak m_1(\mathfrak m(x, y),z)+\mathfrak m(\mathfrak m_1(x, y),z)-\mathfrak m_1(x,\mathfrak m(y,z))-\mathfrak m(x,\mathfrak m_1(y,z))\\
-&\mathfrak m_1(\mathfrak m(y, x),z)-\mathfrak m(\mathfrak m_1(y, x),z)+\mathfrak m_1(y,\mathfrak m(x, z))+\mathfrak m(y,\mathfrak m_1(x, z))=0.
\end{align*}
Or equivalently $\delta(\mathfrak m_1)=[\mf,\mf_1]_{MN}=0$.

The $1$-degree multiderivation $\mathfrak m_1$ is called the infinitesimal of the deformation $\mathfrak{m}_t$. More generally, if $\mathfrak m_i=0$ for $1\leq i\leq n-1$ and $\mathfrak m_n$
 is non zero $1$-degree multiderivation then $\mathfrak m_n$ is called the $n$-infinitesimal of the deformation $\mathfrak{m}_t$.
By the above discussion, the following proposition follows immediately.

\begin{proposition}
The infinitesimal of the deformation $\mathfrak{m}_t$ is a 2-cocycle in $\mathfrak D^1(L)$. More generally, the $n$-infinitesimal is a 2-cocycle.
\end{proposition}
Now we give a notion of equivalence of two deformations.
Let us denote a deformation $(L,A,\mf_t, \ell_t)$ of  $(L,A,\mf, \ell)$
simply by  $L_t$. Let us consider two deformations $L_t$ and $L'_t$ 
of $(L,A,\mf, \ell)$, generated by $\mf_i$ and $\mf'_i$, respectively, for $i\geq 0$.

\begin{definition}
Two deformations $L_t$ and $L'_t$ are said to be equivalent if there exists a formal automorphism $$\Phi_t : L[[t]] \rightarrow L[[t]] \text{ defined as }
\Phi_t=id_L+\sum_{i\geq 1}t^i\phi_i$$
where for each $i\geq 1$, $\phi_i:L \rightarrow L$ is a $\mathbb K$-linear map such that
$$
 {\mathfrak{m}}'_t (x,y)=\Phi_t^{-1}\mathfrak{m}_t(\Phi_t(x),\Phi_t(y))
~~~~\hspace{0.5cm}\text{and}~~~~\hspace{0.5cm}\ell'_t(\Phi_t(x))=\ell_t(x).$$
\end{definition}
\begin{definition}
Any deformation that is equivalent to the deformation $\mathfrak  m_0=\mathfrak{m}$ is said to be a trivial deformation.
\end{definition}
\begin{theorem}\label{Eq1}
The cohomology class of the infinitesimal of a deformation $\mathfrak{m}_t$ is determined by the equivalence class of $\mathfrak{m}_t$.
\end{theorem}

\begin{proof}
Let $\Phi_t$ be an equivalence of deformation between $\mathfrak{m}_t $ and $ \tilde{\mathfrak{m}}_t$. Then we get,
$$
\tilde{\mathfrak{m}}_t (x,y)=\Phi_t^{-1}\mathfrak{m}_t(\Phi_tx,\Phi_ty).
$$
Comparing the coefficients of $t$ from both sides of the above equation we have
\begin{align*}
&\tilde{\mathfrak{m}}_1 (x,y)+\Phi_1(\mathfrak{m}_0(x,y))=\mathfrak{m}_1(x,y)+\mathfrak{m}_0(\Phi_1(x),y)+\mathfrak{m}_0(x,\Phi_1(y)),
\end{align*}
or equivalently,
$$
\mathfrak{m}_1-\tilde{\mathfrak{m}}_1=\delta(\phi_1).
$$
This establishes the result.
\end{proof}
\begin{definition}
A left-symmetric Rinehart algebra is said to be rigid if and only if every deformation of it is trivial.
\end{definition}

\begin{theorem}
A non-trivial deformation of a  left-symmetric Rinehart algebra is equivalent to a deformation whose n-infinitesimal is not a coboundary for some $n\geq 1$.
\end{theorem}

\begin{proof}
Let $\mathfrak{m}_t$ be a deformation of left-symmetric Rinehart algebra with n-infinitesimal $ \mathfrak m_n$ for some $n\geq 1$.
Assume that there exists a 2-cochain $\phi \in C^1(L,L)$ with $\delta(\phi)=\mathfrak m_n$. Then set
$$
\Phi_t=id_{L}+\phi t^n ~~\mbox{and define}~~\bar{\mathfrak{m}_t}=\Phi_t \circ \mathfrak{m}_t \circ \Phi_t^{-1}.
$$
Then by computing the expression and comparing coefficients of $t^n$, we get
$$
\bar{\mathfrak{m}}_n-\mathfrak  m_n=-\delta(\phi).
$$
So, ${\bar{\mathfrak{m}}_n}=0$.
We can repeat the argument to kill off any infinitesimal, which is a coboundary.
\end{proof}

\begin{corollary}
If $H^2(L,L)=0$, then all deformations of $L$ are equivalent to a trivial deformation.
\end{corollary}

\subsection{Obstructions to the extension theory of deformations}
Let $(L,A,\cdot,\ell)$ be a  left-symmetric Rinehart algebra. Now we consider the problem of extending a deformation of $\mathfrak{m}$ of
order $n$ to a deformation of $\mathfrak{m}$ of order $(n+1)$. Let $\mathfrak m_t$ and $\ell_t$ be a deformation of order $n$ of $\mathfrak{m}$ and $\ell$ respectively.
That is  $${\mathfrak m}_t=\sum^n_{i=0} \mathfrak m_i t^i=\mathfrak m+\sum^n_{i=1} \mathfrak m_i t^i\quad \text{and}\quad {\ell}_t=\sum^n_{i=0} \mathfrak \ell_i t^i=\ell+\sum^n_{i=1} \ell_i t^i,$$
where $\mathfrak m_i\in \mathfrak{D}^1(L)$ and $\ell_i:L\to Der(A)$ a linear map for each $1\leq i\leq n$ such that
\begin{align}
   &\mathfrak m_i(\mathfrak m_j(x,y),z)-\mathfrak m_i(x,\mathfrak m_j(y,z))=\mathfrak m_i(\mathfrak m_j(y, x),z)-\mathfrak m_i(y,\mathfrak m_j(x, z)),\label{Defe} \\
   & \mf_i(ax,y)=a\mf_i(x,y),\label{Def1}\\
&\mf_i(x,ay)=a\mf_i(x,y)+\ell_i(x)ay\label{Def3}
\end{align}
 for all $1\leq i,j\leq n$. If there exists a 2-cochain $\mathfrak m_{n+1}\in \mathfrak{D}^1(L)$ and $\ell_{n+1}:L\to Der(A)$ such that $(L,A,\tilde{\mathfrak{m}_t},\tilde{\ell_t})$ is a deformation
 of $(L,A,\mathfrak{m},\ell)$ of order $n+1$, where $$\tilde{\mathfrak{m}_t}=\mathfrak{m}_t+\mathfrak m_{n+1} t^{n+1}\quad \text{and} \quad\tilde{\ell_t}=\ell_t+\ell_{n+1} t^{n+1}.$$ Then we say that $\mathfrak{m}_t$ extends to a deformation of order $(n+1)$.
In this case  $\mathfrak{m}_t$ is called extendable.

\begin{definition}
Let $\mathfrak{m}_t$ be a deformation of $\mathfrak{m}$ of order $n$. Consider the cochain
in $C^3(L,L)$ defined as
{\small\begin{align}\label{Obs}
\begin{split}
    Obs_L(x,y,z)
=\sum_{\substack{i+j=n+1;\\ i,j>0}} \Big( \mathfrak m_i(\mathfrak m_j(x,y),z)-\mathfrak m_i&(x,\mathfrak m_j(y,z))  \\
&-\mathfrak m_i(\mathfrak m_j(y, x),z)+\mathfrak m_i(y,\mathfrak m_j(x, z)))\Big),
\end{split}
\end{align}}
for $x,y,z \in L$. The 3-cochain $Obs_L$ is called an obstruction cochain for extending the deformation of $\mathfrak{m}$ of order
 $n$ to a deformation of order $n+1$.
\end{definition}
A straightforward computation gives the following
\begin{proposition}
The obstructions are left-symmetric Rinehart algebra 3-cocycles.
\end{proposition}
\begin{theorem}\label{Obst}
Let $\mathfrak{m}_t$ be a deformation of $\mathfrak{m}$ of order $n$. Then $\mathfrak{m}_t$ extends to a deformation of order $n+1$ if and only
 if the cohomology class of $Obs_L$ vanishes.
\end{theorem}
\begin{proof}
Suppose that a deformation $\mathfrak{m}_t$ of order $n$ extends to a deformation of order $n+1$. Then
$$
\sum_{\substack{i+j=n+1;\\ i,j\geq0}} \Big( \mathfrak m_i(\mathfrak m_j(x,y),z)-\mathfrak m_i(x,\mathfrak m_j(y,z))
-\mathfrak m_i(\mathfrak m_j(y, x),z)+\mathfrak m_i(y,\mathfrak m_j(x, z)))\Big)=0.
$$
As a result, we get $Obs_L=\delta(m_{n+1})$. So,  the cohomology class of $Obs_L$ vanishes.

Conversely, let $Obs_L$ be a coboundary. Suppose that
$$
Obs_L=\delta(\mathfrak m_{n+1})
$$
for some 2-cochain $\mathfrak m_{n+1}$. Define a map $\tilde{\mathfrak{m}_t}:L[[t]]\times L[[t]]\rightarrow L[[t]]$ as follows
$$
\tilde{\mathfrak{m}_t}=\mathfrak{m}_t+\mathfrak m_{n+1}t^{n+1}.
$$
Then for any $x,y,z\in L $, the map $\tilde{\mathfrak{m}_t}$ satisfies the following identity
$$
\sum_{\substack{i+j=n+1;\\ i,j\geq0}} \Big( \mathfrak m_i(\mathfrak m_j(x,y),z)-\mathfrak m_i(x,\mathfrak m_j(y,z))
-\mathfrak m_i(\mathfrak m_j(y, x),z)+\mathfrak m_i(y,\mathfrak m_j(x, z)))\Big)=0.
$$
 This, in turn, implies that $\tilde{\mathfrak{m}_t}$ is a deformation of $\mathfrak{m}$ extending $\mathfrak{m}_t$.
\end{proof}

\begin{corollary}
If $H^3(L,L)=0,$ then every 2-cocycle in $C^2(L,L)$ is the infinitesimal of some deformation of $\mathfrak{m}$.
\end{corollary}

\subsection{Trivial deformation}

 We study deformations of
left-symmetric Rinehart algebras using the deformation cohomology.
Let $(L,A,\cdot,\ell)$ be a left-symmetric Rinehart algebra, and
$\mathfrak m\in C^2(L,L)$. Consider a $t$-parameterized family
of multiplications $\mathfrak m_t:L[[t]] \otimes L[[t]]\rightarrow L[[t]]$
and linear maps $\ell_t: L \to Der(A)$
given by
\begin{align}
 \mathfrak m_t(x,y)=& x\cdot y+t\mathfrak m(x,y), \\
\ell_t=& \ell +  t \Xi_{\mf}.
\end{align}
If  $L_t=(L,A,\mathfrak m_t,\ell_t)$ is a left-symmetric Rinehart algebra for all $t$, we say that $\mathfrak m$ generates a $1$-parameter infinitesimal
deformation of $(L,A,\cdot,\ell)$

Since $\mathfrak m$ is a $2$-cochain, we have
$$\mathfrak m(ax,y)=a\mathfrak m(x,y),\ \text{and}\ 
\mathfrak m(x,ay)=a\mf_t(x,y)+\Xi_{\mf}(x)(a)y,$$ 
which implies that conditions \eqref{compatibility-left-symmetric-Rinehart1} and \eqref{compatibility-left-symmetric-Rinehart2} in Definition
 \ref{defi:left-symmetric-Rinehart} are satisfied for  $\mathfrak m_t$. Then we can deduce that  $(L,A,\mathfrak m_t,\ell_t)$ is
a deformation of $(L,A,\cdot,\ell)$  if and only if
  \begin{eqnarray}\label{2-closed}
x\cdot\mathfrak m(y,z)-y\cdot\mathfrak m(x,z)+\mathfrak m(y,x)\cdot z-\mathfrak m(x,y)\cdot z    \nonumber \\
=\mathfrak m(y,x\cdot z)-\mathfrak m(x,y\cdot z)-\mathfrak m([x,y],z),
\end{eqnarray}
and
\begin{eqnarray}
\mathfrak m(\mathfrak m(x,y),z)-\mathfrak m(x,\mathfrak m(y,z))&=&\mathfrak m(\mathfrak m(y,x),z)-\mathfrak m(y,\mathfrak m(x,z))\label{omega bracket}.
\end{eqnarray}
Equation $(\ref{2-closed})$ means that $\mathfrak m$ is a $2$-cocycle, and equation $(\ref{omega bracket})$ means that $(L,A,\mathfrak m,\Xi_{\mf})$ is a left-symmetric Rinehart algebra.

Recall that a deformation is said to be
 trivial  if there exists a family of left-symmetric Rinehart algebra isomorphisms
$\Id+tN:L_t\longrightarrow L$.

By direct computations, $L_t$ is trivial  if and only if
\begin{eqnarray}
\mathfrak m(x,y)&=&x\cdot N(y)+N(x)\cdot y-N(x\cdot y),\label{Nij1}\\
N\mathfrak m(x,y)&=&N(x)\cdot N(y),\label{Nij2}\\
\ell \circ N &=& \Xi_{\mf}.\label{Nij3}
\end{eqnarray}
Again, equation \eqref{Nij3}  can be obtained from equation \eqref{Nij1}. 
 It follows from $(\ref{Nij1})$ and $(\ref{Nij2})$ that $N$ must
satisfy the following integrability condition
\begin{eqnarray}
N(x)\cdot N(y)-x\cdot N(y)-N(x)\cdot y+N^2(x\cdot
y)=0.\label{integral condition of Nij}
\end{eqnarray}

Now we give the following definition.
\begin{definition}
An  $A$-linear map
$N:L\longrightarrow L$ is called a   Nijenhuis operator on a left-symmetric Rinehart algebra $(L,A,\cdot,\ell)$  if  the Nijenhuis condition
 $(\ref{integral condition of Nij})$ holds.
\end{definition}
Obviously, any Nijenhuis operator on a left-symmetric Rinehart algebra is also a
Nijenhuis operator on the corresponding  sub-adjacent Lie-Rinehart algebra.

We have seen that a trivial deformation of a left-symmetric Rinehart algebra gives rise to a Nijenhuis operator. In fact, the converse is also true as can be seen from the following theorem.

\begin{theorem}\label{theorem:trivial def}
  Let $(L,A,\cdot,\ell)$ be a left-symmetric Rinehart algebra and $N$ be a Nijenhuis operator. Then a
  deformation of  $(L,A,\cdot,\ell)$ can be   obtained by putting
  $$
\mathfrak m(x,y)= \delta N(x,y).
  $$
  Furthermore, this deformation is trivial.
\end{theorem}
\begin{proof} Since $\mathfrak m$ is a coboundary, then it is a cocycle, i.e. equation
\eqref{2-closed} holds. To see that $\mathfrak m$ generates a
deformation, we only need to show that \eqref{omega bracket} holds,
which follows from the Nijenhuis condition \eqref{integral condition
of Nij}. At the end, we can easily check that 
$$
(\Id +tN)(x\cdot_t y)=(\Id+tN) (x)\cdot (\Id+tN)(y),\quad \ell\circ
(\Id+tN)=\ell_t,
$$
which implies that the deformation is trivial. \end{proof}

\begin{theorem}
  Let $(L,A,\cdot,\ell)$ be a left-symmetric Rinehart algebra and $N$ be a Nijenhuis operator. Then $(L,A,\cdot_N,\ell_N=\ell\circ N)$ is a left-symmetric
  Rinehart algebra, where
 $$x\cdot_N y=x\cdot N(y)+N(x)\cdot y-N(x\cdot y),\forall x,y\in L.$$
\end{theorem}
\begin{proof}
It is obvious to show that $(L,\cdot_N)$  is a left-symmetric algebra and  $\ell_N$ is a representation of L on $Der(A)$. Evidently, we have
$$
\ell_N(ax)=a\ell_N(x),\forall x\in L, a\in A.
$$

Furthermore,  for any $x,y\in L$ and $a\in A$ we have
\begin{align*}
x\cdot_N (ay)=&x\cdot N(ay)+N(x)\cdot (ay)-N(x\cdot (ay))\\
=&a(x\cdot N(y))+\ell(x)aN(y)+ a(N(x)\cdot y)+\ell (N(x))ay-aN(x\cdot y)-N(\ell(x)ay)\\
=&a(x\cdot N(y)+N(x)\cdot y-N(x\cdot y))+\ell_N(x)ay+N(\ell(x)ay)-N(\ell(x)ay).\\
=&a(x\cdot_N y)+\ell_N(x)ay.
\end{align*}
Moreover,
\begin{align*}
(ax)\cdot_N y=&(ax)\cdot N(y)+N(ax)\cdot y-N((ax)\cdot y)\\
=&a(x\cdot N(y))+a(N(x)\cdot y)-aN(x\cdot y)\\
=&a(x\cdot N(y)+N(x)\cdot y-N(x\cdot y))\\
=&a(x\cdot_N y).
\end{align*}
Then, $(L,A,\cdot_N,\ell_N=\ell\circ N)$ is a a left-symmetric
  Rinehart algebra.\end{proof}

  Immediately, we have the following result.
\begin{lemma}\label{lem:Nij1}
 Let $(L,A,\cdot,\ell)$ be a left-symmetric Rinehart algebra and $N$ be a Nijenhuis operator. Then for arbitrary positive $j,k\in\Nat$,  the following equation holds
\begin{equation}
N^j(x)\cdot N^k(y)-N^k(N^j(x)\cdot y)-N^j(x\cdot N^k(y))+N^{j+k}(x\cdot y)=0,\quad\forall~ x,y\in L.
\end{equation}
If $N$ is invertible, this formula becomes valid for arbitrary $j,k\in\mathbb Z$.
\end{lemma}

By direct calculations, we have the following corollary.
\begin{corollary}\label{lem:Niejproperty}
  Let $(L,A,\cdot,\ell)$ be a left-symmetric Rinehart algebra and $N$ a Nijenhuis operator.
  \begin{itemize}
\item[$\rm(i)$] For all $k\in\Nat$, $(L,A,\cdot_{N^k},\ell_{N^k}=\ell\circ N^k)$ is a left-symmetric Rinehart algebra.
\item[$\rm(ii)$]For all $l\in\Nat$, $N^l$ is  a Nijenhuis operator on the left-symmetric Rinehart algebra $(L,A,\cdot_{N^k},\ell_{N^k})$.
\item[$\rm(iii)$]The left-symmetric Rinehart algebras $(L,A,(\cdot_{N^k})_{N^l},\ell_{N^{k+l}})$ and $(L,A,\cdot_{N^{k+l}},\ell_{N^{k+l}})$ are the same.
\item[$\rm(iv)$]$N^l$ is a left-symmetric Rinehart algebra homomorphism from $(L,A,\cdot_{N^{k+l}},\ell_{N^{k+l}})$ to $(L,A,\cdot_{N^k},\ell_{N^{k}})$.
  \end{itemize}
\end{corollary}

\begin{theorem}
 Let $(L,A,\cdot,\ell)$ be a left-symmetric Rinehart algebra and $N$ be a Nijenhuis operator. Then the operator $P(N)=\sum_{i=0}^nc_iN^i$ 
 is a Nijenhuis
operator. If $N$ is invertible,  then $Q(N)=\sum_{i=-m}^nc_iN^i$ is also a Nijenhuis
operator.
\end{theorem}
\begin{proof} According to  Lemma \ref{lem:Nij1}, we obtain, $\forall x,y\in L,$
\begin{eqnarray*}
&&P(N)(x)\cdot  P(N)(y)-P(N)(P(N)(x)\cdot  y)-P(N)(x\cdot  P(N)(y))+P^2(N)(x\cdot  y)\\
&=&\sum_{i,j=0}^nc_jc_k\Big(N^j(x)\cdot  N^k(y)-N^k(N^j(x)\cdot  y)-N^j(x\cdot  N^k(y))+N^{j+k}(x\cdot  y)\Big)=0,  
\end{eqnarray*}

which implies that $P(N)$ is a Nijenhuis operator. Similarly we can easy check the second statement.\end{proof}

\section{\texorpdfstring{$\huaO$}{O}-operators and Nijenhuis
operators}

In this section, we highlight  the relationships between
$\huaO$-operators and Nijenhuis operators on left-symmetric Rinehart algebras.  Moreover, we illustrate some connections
  between Nijenhuis operators and  compatible $\huaO$-operators  on
left-symmetric Rinehart algebras.

\subsection{Relationships between \texorpdfstring{$\huaO$}{O}-operators and Nijenhuis
operators}

We first give the definitions of an {$\huaO$-operator} and of Rota-Baxter operator.
\begin{definition}
An {$\huaO$-operator} on a left-symmetric Rinehart algebra
$(L,A,\cdot,\ell)$ associated to a representation $(M;\rho,\mu)$ is a
linear map $T:M\longrightarrow L$ satisfying
\begin{eqnarray}
  T(au)&=& aT(u),\\ \label{o-operator1}
  T(u)\cdot  T(v)&=&T\Big(\rho(T(u))(v)+\mu(T(v))(u)\Big),\quad \forall u,v\in M, a\in A.\label{o-operator2}
\end{eqnarray}
\end{definition}

\begin{definition}
Let $(L,A,\cdot,\ell)$ be a left-symmetric Rinehart algebra and $\huaR:ker(\ell)\longrightarrow
L$ a   linear operator. If $\huaR$ satisfies
\begin{eqnarray}
  \huaR(ax)&=& a\huaR(x),\label{rota-baxter1}\\
\huaR(x)\cdot  \huaR(y)&=&\huaR(\huaR(x)\cdot  y+x\cdot  \huaR(y)),\quad \forall x,y\in L, a\in A,\label{rota-baxter2}
\end{eqnarray}
 then $\huaR$ is called a
 Rota-Baxter operator of weight $0$ on $L$.
\end{definition}
  Notice that a Rota-Baxter operator of weight zero on a
left-symmetric Rinehart algebra $L$ is exactly an $\huaO$-operator
associated to the adjoint representation $(L;ad^L,ad^R)$.

The following proposition gives connections between Nijenhuis operators and Rota-Baxter operators.

\begin{proposition}\label{pro:RN}
Let $(L,A,\cdot,\ell)$ be a  left-symmetric Rinehart algebra and $N:L\longrightarrow
L$ a linear operator.
\begin{enumerate}
\item[\rm(i)] If $N^2={\rm Id}$, then $N$ is a
Nijenhuis operator if and only if $N\pm{\rm Id}$ is a Rota-Baxter
operator of weight $\mp 2$ on $(L,A,\cdot,\ell)$.
\item[\rm(ii)] If $N^2=0$, then
$N$ is a Nijenhuis operator if and only if $N$ is a Rota-Baxter
operator of weight zero on $(L,A,\cdot,\ell)$.
\item[\rm(iii)] If $N^2=N$, then $N$ is a
Nijenhuis operator if and only if $N$ is a Rota-Baxter operator of
weight $-1$ on $(L,A,\cdot,\ell)$.
\end{enumerate}

\end{proposition}

\begin{proof}
 For Item (i), for all $x,y\in L$ then we have
 \begin{align*}
 (N-Id)&(x)\cdot (N-Id)(y)\\
 &\qquad-(N-Id)\big((N-Id)(x)\cdot y+x\cdot (N-Id)(y)\big)+2 (N-Id)(x\cdot y)\\
 &=N(x)\cdot N(y)-N\big(N(x)\cdot y+x\cdot N(y)\big)+x\cdot y \\
  &=N(x)\cdot N(y)-N\big(N(x)\cdot y+x\cdot N(y)\big)+N^2(x\cdot y).
 \end{align*}
 So  $N$ is a Nijenhuis operator if and only if $N- Id$ is a Rota-Baxter
operator of weight $2$ on $L$.
Similarly, we obtain that $N$ is a Nijenhuis operator if and only if $N+ Id$ is a Rota-Baxter
operator of weight $-2$ on $L$.

Items (ii) and (iii) are obvious from the definitions of Nijenhuis operators and  Rota-Baxter
operator.
\end{proof}

\begin{proposition}\label{lem:OR}  Let $(L,A,\cdot,\ell)$ be a  left-symmetric Rinehart algebra  and
$(M;\rho,\mu)$ be a representation on $L$. Let $T:M\rightarrow L$ be a
linear map. For any $\lambda$, $T$ is an $\huaO$-operator on
$L$ associated to $(M;\rho,\mu)$ if and only if the linear map
$\huaR_{T,\lambda}:=\left(\begin{array}{cc}0&T\\0&-\lambda {\rm
Id}\end{array}\right)$ is a Rota-Baxter operator of weight
$\lambda$ on the semidirect product left-symmetric Rinehart algebra
$({L\oplus M},\cdot_{L\oplus M})$, where the multiplication $\cdot_{L\oplus M}$
is given by \eqref{eq:mulsemi}.
\end{proposition}
\begin{proof}
It is easy to check the equation \eqref{rota-baxter1}. Let $x_1,x_2\in L$ and $m_1,m_2\in M$,
\begin{align}
\huaR_{T,\lambda}(x_1+m_1)\cdot_{L\oplus M}  \huaR_{T,\lambda}(x_2+m_2)&= (T(m_1)-\lambda m_1)\cdot_{L\oplus M} (T(m_2)-\lambda m_2)\nonumber\\
&=T(m_1)\cdot T(m_2)-\lambda\rho(T(m_1)m_2-\lambda\mu(T(m_2))m_1.\label{eq1}
\end{align}
On the other hand,
{\allowdisplaybreaks
\begin{align}
&\huaR_{T,\lambda}\big(\huaR_{T,\lambda}(x_1+m_1)\cdot_{L\oplus M}  (x_2+m_2)+(x_1+m_1)\cdot_{L\oplus M}  \huaR_{T,\lambda}(x_2+m_2)\big)+ \nonumber\\
&
\quad \lambda
\huaR_{T,\lambda}((x_1+m_1)\cdot_{L\oplus M}  (x_2+m_2))\nonumber\\
=&\huaR_{T,\lambda}\big( (T(m_1)-\lambda m_1)\cdot_{L\oplus M}(x_2+m_2) + (x_1+m_1)\cdot_{L\oplus M}(T(m_2)-\lambda m_2)\big)+ \nonumber\\
&
 \quad \lambda
\huaR_{T,\lambda}(x_1\cdot x_2+\rho(x_1)m_2+\mu(x_2)m_1)\nonumber\\ \pagebreak
=&\huaR_{T,\lambda}\big( T(m_1)\cdot x_2+\rho(T(m_1))m_2-\lambda\mu(x_2)m_1\\ 
& \quad\qquad+x_1\cdot T(m_2)-\lambda\rho(x_1)m_2+\mu(T(m_2))m_1\big)\nonumber\\
&+\lambda \big(T(\rho(x_1)m_2)+T(\mu(x_2)m_1)-\lambda(\rho(x_1)m_2+\mu(x_2)m_1)\big)\nonumber\\
\pagebreak
=& T\big(\rho(T(m_1))m_2+\mu(T(m_2))m_1\big)-\lambda T\big(\mu(x_2)m_1+\rho(x_1)m_2\big)\\
&-\lambda \Big( \rho(T(m_1))m_2+\mu(T(m_2))m_1\Big)+\lambda^2 \Big(\mu(x_2)m_1+\rho(x_1)m_2\big)\nonumber\\
&+\lambda \big(T(\rho(x_1)m_2)+T(\mu(x_2)m_1)\big)-\lambda^2\big(\rho(x_1)m_2+\mu(x_2)m_1)\big)\nonumber\\
=& T\big(\rho(T(m_1))m_2+\mu(T(m_2))m_1\big)-\lambda  \rho(T(m_1))m_2-\lambda\mu(T(m_2))m_1.\label{eq2}
\end{align}}

According to equations \eqref{eq1} and \eqref{eq2},  $\huaR_{T,\lambda}$ is a Rota-Baxter operator of weight
$\lambda$ on the semidirect product left-symmetric Rinehart algebra
$({L\oplus M},\cdot_{L\oplus M})$ if and only if  $T$ is an $\huaO$-operator on
$(L,A,\cdot,\ell)$ associated to $(M;\rho,\mu)$.
\end{proof}
\begin{proposition}
Let $(L,A,\cdot,\ell)$ be a  left-symmetric Rinehart algebra and let $(M;\rho,\mu)$ be a
representation on $L$. Let $T:M\rightarrow L$ be a linear map. Then the
following statements are equivalent.
\begin{enumerate}
\item[\rm(i)] $T$ is an $\huaO$-operator on the left-symmetric Rinehart algebra $(L,A,\cdot,\ell)$.
\item[\rm(ii)]
  $\huaN_T :=\left(\begin{array}{cc}0&T\\0&{\rm Id}\end{array}\right)$ is a Nijenhuis operator on the
left-symmetric Rinehart algebra $({L\oplus M},\cdot_{L\oplus M})$.
\item[\rm(iii)]
  $N_T :=\left(\begin{array}{cc}0&T\\0&0\end{array}\right)$ is a Nijenhuis operator on
the left-symmetric Rinehart algebra  $({L\oplus M},\cdot_{L\oplus M})$.
\end{enumerate}

\end{proposition}

\begin{proof} Note that $\huaN_T=\huaR_{T,-1}$  and
$(\huaN_T )^2=\huaN_T$, thus $\huaN_T$ is a Nijenhuis operator on
the left-symmetric Rinehart algebra  $({L\oplus M},\cdot_{L\oplus M})$, using (iii) in Proposition \ref{pro:RN}.

Similarly $N_T=\huaR_{T,0}$ and
$(N_T)^2=0$, then $N_T$ is a Nijenhuis operator on
the left-symmetric Rinehart algebra  $({L\oplus M},\cdot_{L\oplus M})$, according to (ii) in Proposition \ref{pro:RN}. \end{proof}

\subsection{Compatible \texorpdfstring{$\huaO$}{O}-operators and Nijenhuis operators}
In this subsection we study compatibility of $\huaO$-operators and Nijenhuis operators.  First we start with the following definition.

\begin{definition}
Let $(L,A,\cdot,\ell)$ be a  left-symmetric Rinehart algebra and let $(M;\rho,\mu)$ be a
representation. Let $T_1,T_2: M\longrightarrow L$ be two
$\huaO$-operators associated to $(M;\rho,\mu)$. Then $T_1$ and $T_2$ are called
compatible if  $T_1+T_2$ is an $\huaO$-operator
associated to $(M;\rho,\mu)$.
\end{definition}

Let $T_1,T_2: M\longrightarrow L$ be two $\huaO$-operators
on a left-symmetric Rinehart algebra $(L,A,\cdot,\ell)$ associated to a
representation $(M;\rho,\mu)$ such that
{\begin{align}
 T_1(u)\cdot  T_2(v)+ T_2(u)\cdot
 T_1(v)&=T_1\Big(\rho(T_2(u))(v)+\mu(T_2(v))(u)\Big)  \nonumber \\
 &\qquad+T_2\Big(\rho(T_1(u))(v) +\mu(T_1(v))(u)\Big),\label{eq:CN1}
\end{align}}
for all $u,v\in M.$

\begin{lemma}Two operators $T_1$ and $T_2$ are compatible
if and only if the  equation \eqref{eq:CN1} holds.
\end{lemma}

\begin{proof} For all $u,v\in M$, $a\in A$, we have
\begin{align*}
 (T_1+T_2)(au)&=T_1(au)+T_2(au)\\
 &=aT_1(u)+aT_2(u)\\
 &=a(T_1+T_2)(u).  
\end{align*}
Furthermore,
\begin{align*}
 &\!\!(T_1+T_2)(u)\cdot (T_1+T_2)(v)-(T_1+T_2)\Big(\rho((T_1+T_2)(u))(v)+\mu((T_1+T_2)(v))(u)\Big)\\
 &=T_1(u)\cdot T_1(v)+T_1(u)\cdot T_2(v)+T_2(u)\cdot T_1(v)+T_2(u)\cdot T_2(v)\\
 &\qquad-(T_1+T_2)\Big(\rho(T_1(u))(v)+\rho(T_2(u))(v)+\mu(T_1(v))(u)+\mu(T_2(v))(u)\Big)\\
 &=T_1(u)\cdot T_1(v)+T_1(u)\cdot T_2(v)+T_2(u)\cdot T_1(v)+T_2(u)\cdot T_2(v)\\
 &\qquad-T_1\big(\rho(T_1(u))(v)+\mu(T_1(v))(u)\big)-T_1\big(\rho(T_2(u))(v)+\mu(T_2(v))(u)\big)\\
 &\qquad-T_2\big(\rho(T_1(u))(v)+\mu(T_1(v))(u)\big)-T_2\big(\rho(T_2(u))(v)+\mu(T_2(v))(u)\big)\\
 &=T_1(u)\cdot  T_2(v)+ T_2(u)\cdot T_1(v) \\
&\qquad-T_1\Big(\rho(T_2(u))(v)+\mu(T_2(v))(u)\Big) -T_2\Big(\rho(T_1(u))(v)+\mu(T_1(v))(u)\Big)
\end{align*}
Then $T_1+T_2$ is an $\huaO$-operator associated to $(M;\rho,\mu)$ if and only if equation \eqref{eq:CN1} holds.
\end{proof}
\begin{remark}
Equation \eqref{eq:CN1} implies that for any $k_1, k_2$  the linear combination $k_1 T_1+k_2 T_2$  is an $\huaO$-operator.
\end{remark}
There is a close relationship between a Nijenhuis operator and a pair of compatible
$\huaO$-operators as can be seen from the following proposition.
\begin{proposition}\label{pro:TTN}
Let $T_1,T_2: M\longrightarrow L$ be two $\huaO$-operators on
a left-symmetric Rinehart algebra $(L,A,\cdot,\ell)$ associated to a representation
$(M;\rho,\mu)$. Suppose that $T_2$ is invertible. If $T_1$ and $T_2$
are compatible, then $N=T_1\circ T_2^{-1}$ is a Nijenhuis operator on the left-symmetric Rinehart algebra $(L,A,\cdot,\ell)$.
\end{proposition}

\begin{proof} For all $x,y\in L$, since $T_2$ is invertible, there exist $u,v\in M$ such that
$T_2(u)=x,$ $T_2(v)=y$. Hence $N=T_1\circ T_2^{-1}$ is a Nijenhuis operator
if and only if the following equation holds:
\begin{equation}\label{eq11}
NT_2(u)\cdot  NT_2(v)=N(NT_2(u)\cdot  T_2(v)+T_2(u)\cdot
NT_2(v))-N^2(T_2(u)\cdot  T_2(v)).
\end{equation}
Since $T_1=N\circ T_2$ is an
$\huaO$-operator, the left-hand side of the above equation is
$$NT_2(\rho(NT_2(u))(v)+\mu(NT_2(v))(u)).$$
Using the fact  that $T_2$  and $T_1=N\circ T_2$ are two compatible
$\huaO$-operators, we get
\begin{eqnarray*}
 &&NT_2(u)\cdot  T_2(v)+T_2(u)\cdot  NT_2(v)\\
 &=&
T_2(\rho(NT_2(u))(v)+\mu(NT_2(v))(u))+NT_2(\rho(T_2(u))(v)+\mu(T_2(v))(u))\\
&=&T_2(\rho(NT_2(u))(v)+\mu(NT_2(v))(u))+N(T_2(u)\cdot  T_2(v)).
\end{eqnarray*}
Hence, equation \eqref{eq11} holds by acting  $N$ on both sides of the last equality. \end{proof}

Using an $\huaO$-operator and a Nijenhuis operator, we can construct a pair of compatible $\huaO$-operators.

\begin{proposition}\label{pro:NT}
Let $T: M\longrightarrow L$ be an $\huaO$-operator on a left-symmetric Rinehart algebra $(L,A,\cdot,\ell)$ associated to a representation
$(M;\rho,\mu)$ and let $N$ be a Nijenhuis operator on $(L,A,\cdot,\ell)$.
Then $N\circ T$ is an $\huaO$-operator on the left-symmetric Rinehart algebra $(L,A,\cdot,\ell)$ associated to $(M;\rho,\mu)$
if and only if for all $u,v\in M$, the following equation holds:
{\small\begin{align}
&N\Big(NT(u)\cdot  T(v)+T(u)\cdot  NT(v)\Big)\nonumber\\
=&N\Big(T\big(\rho(NT(u))(v)+\mu(NT(v))(u)\big)+NT\big(\rho(T(u))(v)+\mu(T(v))(u)\big)\Big).\label{eq:ON}
\end{align}}
If in addition $N$ is invertible, then $T$ and $N T$
are compatible. More explicitly, for any $\huaO$-operator $T$,
if there exists an invertible Nijenhuis operator $N$ such that $N T$
is also an $\huaO$-operator, then $T$ and $N T$ are compatible.
\end{proposition}

\begin{proof}   Let $u,v\in M$ and $a\in A$, we have
\begin{align*}
N T(au)=N(T(au))=N(aT(u))=aN T(u).
\end{align*}
In addition, since $N$ is a Nijenhuis operator and $T$ is an $\huaO$-operator we have
\begin{align*}
NT(u)\cdot  NT(v)&=N\Big(NT(u)\cdot  T(v)+T(u)\cdot
NT(v)\Big)-N^2(T(u)\cdot  T(v))\\
&=NT\Big(\rho(NT(u))(v)+\mu(NT(v))(u)\Big)
\end{align*}
if and only if  (\ref{eq:ON}) holds.

If $N T$ is an $\huaO$-operator and $N$ is invertible, then we have
\begin{align*}
&T(u)\cdot  T(v)+T(u)\cdot  NT(v) \\
&=T\big(\rho(NT(u))(v)+\mu(NT(v))(u)\big)+NT\big(\rho(T(u))(v)+\mu(T(v))(u)\big),\end{align*}
which is exactly the condition that $N T$ and $T$ are compatible. \end{proof}

The following result is an immediate consequence of the last two propositions.
\begin{corollary}
Let $T_1,T_2: M\longrightarrow L$ be two $\huaO$-operators on
a left-symmetric Rinheart algebra $(L,A,\cdot,\ell )$ associated to a representation
$(M;\rho,\mu)$. Suppose that  $T_1$ and $T_2$ are invertible. Then
$T_1$ and $T_2$ are compatible if and only if $N=T_1\circ T_2^{-1}$ is a
Nijenhuis operator.
\end{corollary}

\subsection*{Acknowledgement}
The authors are thankful to the Deanship of Graduate Studies and Scientific Research at University of Bisha for supporting this work through the Fast-Track Research Support Program.  Mohamed Elhamdadi was partially supported by Simons Foundation collaboration grant 712462.

\EditInfo{December 6, 2023}{February 29, 2024}{Adam Chapman and Ivan Kaygorodov}

\end{document}